\documentclass[preprint,12pt]{elsarticle1}

\usepackage{subfig}
\usepackage{graphicx}
\usepackage{caption,color}   % 24.62
\usepackage{amssymb}
\usepackage{amsthm}
\usepackage{amsmath}
\usepackage{epic}
\usepackage{setspace}
\usepackage{float}
\usepackage{multirow}
\newtheorem{thm}{Theorem}[section]

\newtheorem{lem}[thm]{Lemma}

\newtheorem{defi}[thm]{Definition}

\numberwithin{equation}{section}
\journal{}

\begin{document}
\begin{spacing}{1.15}
\begin{frontmatter}
\title{Estrada index and subgraph centrality of hypergraphs via tensors}

\author{Hong Zhou}
\author{Lizhu Sun}\ead{sunlizhu678876@126.com}
\author{Changjiang Bu}
\address{College of Mathematical Sciences, Harbin Engineering University, Harbin 150001, PR China}

\begin{abstract}
Uniform hypergraphs have a natural one-to-one correspondence to tensors. In this paper, we investigate the Estrada index and subgraph centrality of an $m$-uniform hypergraph $\mathcal{H}$ via the adjacency tensor. We establish some bounds for the Estrada index and give expressions of the subgraph centrality in terms of graph parameters of the multi-digraphs associated with $\mathcal{H}$. When $\mathcal{H}$ is $2$-uniform, the above Estrada index and subgraph centrality are the Estrada index and subgraph centrality of a graph.

\end{abstract}

\end{frontmatter}

\section{Introduction}
For a simple undirected graph $H$, its \textit{Estrada index}
\begin{equation*}\label{S2}
EE(H)=\sum \limits_{i=1}^{n}e^{\lambda_{i}},
\end{equation*}
where $\lambda_1, \lambda_{2},\ldots,\lambda_{n}$ are all the eigenvalues of the adjacency matrix of $H$ \cite{ESTRADA2000713}. The study of Estrada index has attracted extensive attention\cite{DELAPENA200770,Gutman2008Lower,2008On,2011The,2012Estrada}. The Estrada index of graphs has wide applications in biology \cite{ESTRADA2000713}, chemistry \cite{molecules} and complex networks \cite{PhysRevE.71.056103,PhysRevE.72.046105}.
The Estrada index of graphs is closely related to the subgraph  centrality of a vertex in graphs and the trace of the adjacency matrix of graphs \cite{DELAPENA200770,PhysRevE.71.056103}.

%\cite{6,21,7,7+,8-}, physics \cite{4}
Let $A$ be the adjacency matrix of the graph $H$ and $\mu_{d}(j)=(A^{d})_{jj}$. Then $\mu_{d}(j)$ is the number of closed walks
of length $d$ starting and ending at the vertex $j$ in $H$ \cite{1980Spectra}.
\begin{align}\label{wa1}
C(j)=\sum \limits_{d=0}^{\infty}\frac{\mu_{d}(j)}{d!}
\end{align}
is called the \textit{subgraph centrality} of a vertex $j$ in $H$ \cite{PhysRevE.71.056103}.
The subgraph centrality is a topological parameter to measure the importance of nodes in networks, which is widely used in the real-world network analysis \cite{PhysRevE.71.056103,2011Thee}.

%which can reflect the centrality of $j$ in $H$ from the perspective of the number of closed walks starting and ending at the vertex $j$ \cite{2}.
The \textit{$d$th order spectral moment} of $H$ is the sum of $d$ powers of all the eigenvalues of $A$,  denoted by $S_{d}(H)$. Since the trace
$tr(A^d)=\sum_{j=1}^n\mu_{d}(j)=S_{d}(H)$ \cite{1980Spectra},
\begin{align}\label{ss1}
\sum_{j=1}^nC(j)=\sum \limits_{d=0}^{\infty}\sum_{j=1}^n\frac{\mu_{d}(j)}{d!}=\sum \limits_{d=0}^{\infty}\frac{tr(A^d)}{d!}=
\sum \limits_{d=0}^{\infty}\sum_{i=1}^n \frac{\lambda_i^d}{d!}=\sum_{i=1}^n e^{\lambda_i}=EE(H).
\end{align}

Since uniform hypergraphs have a natural one-to-one correspondence to tensors, in this paper, we investigate the Estrada index and subgraph centrality of hypergraphs via tensors.
Next, we introduce some notations and concepts for tensors and hypergraphs. Let $[n]=\{1,2,\ldots,n\}$, $[n]^m=\{i_1i_2\cdots i_m|~i_k \in[n],~k=1,\ldots,m\}$
and $\mathbb{C}$ be complex field. An order $m$ dimension $n$ complex tensor

$$
\mathcal{T}=\left( {t_{\alpha} } \right),~\mbox{for }\alpha\in[n]^m, t_{\alpha}\in\mathbb{C},
$$
is a multidimensional array with $n^m$ entries.
%For the convenience of notations we write $t_{i_{1}i_{2}...i_{m}}=t_{i_{1}\alpha}$, where $\alpha=i_{2}\cdot\cdot\cdot i_{m}.$
When $m=2$, $\mathcal{T}$ is an $n\times n$ matrix \cite{qi2005eigenvalues,lim2005singular}.

%\in \mathbb{C}^{[m,n]}
A hypergraph $\mathcal{H}=(V(\mathcal{H}), E(\mathcal{H}))$ is called \textit{$m$-uniform} if $|e|=m\geq2$ for all $e\in E(\mathcal{H})$. For $e=\{i_{1},i_{2},\ldots,i_{m}\}\in E(\mathcal{H})$, $e$ is also written $i_{1}i_{2}\cdots i_{m}$ in this paper. For an $m$-uniform hypergraph $\mathcal{H}$, its  \textit{adjacency tensor} is the order $m$ dimension $n$ tensor
$\mathcal{A}_\mathcal{H}=(h_{\alpha})$, where
%is the \textit{adjacency tensor} of an $m$-uniform hypergraph $\mathcal{H}$ \cite{10}, where
\begin{equation*}
h_{\alpha}=\begin{cases}
\frac{1}{(m-1)!},& \text{if } \alpha\in E(\mathcal{H}),\notag \\
0,& \text{otherwise}.
\end{cases}
\end{equation*}
Clearly, $\mathcal{A}_\mathcal{H}$ is the adjacency matrix of $\mathcal{H}$ when $\mathcal{H}$ is $2$-uniform. The eigenvalues of $\mathcal{A}_\mathcal{H}$ are also called the eigenvalues of $\mathcal{H}$ \cite{cooper2012spectra}.

In 2005, the concept of eigenvalues of tensors was proposed by Qi \cite{qi2005eigenvalues} and Lim \cite{lim2005singular}, independently.
The eigenvalues of tensors and related problems are important research topics of spectral hypergraph theories \cite{cooper2,doi:10.1137/21M1404740,2017Tensor,clark2021harary}, especially the trace of tensors \cite{clark2021harary,2011Analogue,hu2013determinants,shao2015some,doi:10.1080/03081087.2021.1953431}.

%\cite{11,s12+,12,12+,l12+}\cite{1+,17,4+}

%The research of spectral hypergraph theory via tensors has developed rapidly [17-19].

Morozov and Shakirov  gave an expression of the $d$th order trace $Tr_{d}(\mathcal{T})$ of a tensor $\mathcal{T}$ \cite{2011Analogue}. Hu et al. proved that $Tr_{d}(\mathcal{T})$ is equal to the sum of $d$ powers of all eigenvalues of $\mathcal{T}$ \cite{hu2013determinants}. For a uniform hypergraph $\mathcal{H}$, the sum of $d$ powers of all eigenvalues of $\mathcal{A}_\mathcal{H}$ is called the \textit{$d$th order  spectral moment} of $\mathcal{H}$, denoted by $S_{d}(\mathcal{H})$. Then
$Tr_{d}(\mathcal{\mathcal{A}_\mathcal{H}})=S_{d}(\mathcal{H})$. Shao et al. established some formulas for the $d$th order trace of tensors in terms of some graph parameters \cite{shao2015some}.
Clark and Cooper expressed the spectral moments of
hypergraphs by the numbers of Veblen multi-hypergraphs and used this result to give the ``Harary-Sachs'' coefficient theorem for hypergraphs \cite{clark2021harary}.
Chen et al. gave a formula for the spectral moment of a hypertree in terms of the numbers of some subhypertrees \cite{doi:10.1080/03081087.2021.1953431}.

In this paper, we define the Estrada index and subgraph centrality of a uniform hypergraph $\mathcal{H}$ via the adjacency tensor. The bounds for the Estrada index are established.
%The subgraph centrality $C(j)$ of a vertex $j$ in $\mathcal{H}$ is defined by using $\mu_d(j)$, where $\mu_d(j)$ is the term corresponding to the vertex $j$ in $Tr_d(\mathcal{A}_\mathcal{H})$, $d=0,1,2, \ldots .$
%When $\mathcal{H}$ is $2$-uniform, the above subgraph centrality is the subgraph centrality of graphs as in Equation (\ref{wa1}).
We give two expressions of the subgraph centrality by the number of  Eulerian closed walks of the multi-digraphs  associated with $\mathcal{H}$ and the number of arborescences of the multi-digraphs associated with  $\mathcal{H}$, respectively.
Similar to the Estrada index of a graph as in Equation (\ref{ss1}),  the Estrada index of a uniform hypergraph $\mathcal{H}$ is equal to the sum of the subgraph centrality measures of all vertices in $\mathcal{H}$.

\section{Preliminaries}

%In this section, we define the Estrada index of a uniform hypergraph $\mathcal{H}$ and the subgraph centrality of vertices in $\mathcal{H}$. Some concepts and lemmas of necessity are introduced.
%We show that the Estrada index of $\mathcal{H}$ can be obtained by both the sum of the subgraph centrality  measures of all vertices in $\mathcal{H}$ and the sum of all the $d$th order traces $ Tr_d(\mathcal{A}_\mathcal{H})$ divided by $d!$, $d=0,1,2,\ldots.$
Let $\mathbb{C}^{n}$ be the set of $n$-dimension complex vectors and $\mathbb{C}^{[m,n]}$ be the set of complex tensors with order $m$ dimension $n$.
For a tensor $\mathcal{T}=(t_{i\alpha})\in\mathbb{C}^{[m,n]}$ and $x=\left({x_1 ,\ldots ,x_n}\right)^\mathrm{T}\in\mathbb{C}^n$,  $\mathcal{T}x^{m-1}$ is a vector in $\mathbb{C}^n$ whose $i$-th component is
\begin{align*}
(\mathcal{T}x^{m-1})_i=\sum\limits_{\alpha\in[n]^{m-1}}t_{i\alpha}x^{\alpha},
\end{align*}
where $x^{\alpha}=x_{i_{1}}x_{i_{2}}\cdots x_{i_{m-1}}$ for $\alpha=i_{1}i_{2}\cdots i_{m-1}$.
%For a tensor $\mathcal{T}\in\mathbb{C}^{[m,n]}$,
A number $\lambda\in\mathbb{C}$ is called an \textit{eigenvalue} of $\mathcal{T}$ if there exists a nonzero vector $x\in\mathbb{C}^n$ such that
$$\mathcal{T}x^{m-1}=\lambda x^{[m-1]},$$
where $x^{\left[ {m - 1} \right]}  = \left( {x_1^{m - 1} ,\ldots,x_n^{m - 1} } \right)^\mathrm{T}$.
The number of eigenvalues of $\mathcal{T}$ is $k=n(m-1)^{n-1}$ \cite{qi2005eigenvalues,lim2005singular}.
 Let $\rho=\max\{|\lambda_1|,|\lambda_2|,\dots,|\lambda_k| \}$ be the spectral radius of $\mathcal{T}$, where $\lambda_1,\lambda_2,\ldots,\lambda_k$ are all the eigenvalues of $\mathcal{T}$.
%The \textit{spectral radius} of $\mathcal{T}$ is  $\rho=\max\{|\lambda_1|,|\lambda_2|,\dots,|\lambda_k| \}$, where $\lambda_1,\lambda_2,...,\lambda_k$ are all the eigenvalues of $\mathcal{T}$.

%\begin{align}\label{ss2}
%[n]^m=\{i_1i_2\cdots i_m|~i_k \in[n],~k=1,\ldots,m\}.
%\end{align}

The \textit{$d$th order trace} $Tr_{d}(\mathcal{T})$ of a tensor $\mathcal{T}=(t_{\alpha})\in \mathbb{C}^{[m,n]}$ is expressed as follows \cite{2011Analogue}:
\begin{align}\label{waaa2}
&Tr_{d}(\mathcal{T}) \notag\\
&=(m-1)^{n-1}\sum \limits_{d_{1}+\cdots+d_{n}=d}\prod \limits_{i=1}^{n}\frac{1}{(d_{i}(m-1))!}(\sum \limits_{\alpha_{i}\in [n]^{m-1}}t_{i\alpha_{i}}\frac{\partial}{\partial a_{i\alpha_{i}}})^{d_{i}}tr(A^{d(m-1)}),
\end{align}
where $A=(a_{ij})$ is an $n\times n$ auxiliary matrix, $d_{1},\ldots,d_{n}$ are nonnegative integers and
$\frac{\partial}{\partial a_{i\alpha_{i}}}:=\frac{\partial}{\partial a_{ii_{2}}}\cdot \cdot \cdot \frac{\partial}{\partial a_{ii_{m}}}$ for $\alpha_{i}=i_{2}\cdot \cdot \cdot i_{m}$.
%When $\mathcal{T}$ is a matrix, $Tr_{d}(\mathcal{T})=tr(\mathcal{T}^{d})$.
%the above $d$th order trace of $\mathcal{T}$ is the trace of the $d$ power of $\mathcal{T}$    and $[n]^{m-1}$ is defined as in Equation (\ref{ss2}).

In \cite{cooper2012spectra}, the $d$th order traces of the adjacency tensor of an $m$-uniform hypergraph were given for $d=0,1,2,\ldots,m$.
\begin{lem}\cite{cooper2012spectra}\label{L2}
Let $\mathcal{H}$ be an $m$-uniform hypergraph with $n$ vertices and $q$ edges. Then

(1) $Tr_{0}(\mathcal{A}_\mathcal{H})=n(m-1)^{n-1}$;

(2) $Tr_d(\mathcal{A}_\mathcal{H})=0$ for $d=1,2,\ldots, m-1$;

(3) $Tr_{m}(\mathcal{A}_\mathcal{H})=qm^{m-1}(m-1)^{n-m}$.

\end{lem}
%$\mu_{0}(j)=(m-1)^{n-1},j=1,2,...,n.$

Since uniform hypergraphs have a natural one-to-one correspondence to tensors, we define the Estrada index and subgraph centrality of a uniform hypergraph via the adjacency tensor.
\begin{defi}
For an $m$-uniform hypergraph $\mathcal{H}$ with $n$ vertices,
\begin{equation*}
\sum \limits_{i=1}^{k}e^{\lambda_{i}}
\end{equation*}
is called the Estrada index of $\mathcal{H}$, denoted by $EE(\mathcal{H})$, where $\lambda_1,\lambda_{2}, \ldots,\lambda_k$ are all the eigenvalues of $\mathcal{A}_\mathcal{H}$.
\end{defi}

Clearly, when $\mathcal{H}$ is $2$-uniform, the above index is the Estrada index of a graph \cite{ESTRADA2000713}.

Let $\mathcal{H}=(V(\mathcal{H}),E(\mathcal{H}))$ be an $m$-uniform hypergraph with $n$ vertices.
For  $j \in V(\mathcal{H})$, let
$\mu_{d}(j)$ be the term corresponding to the vertex $j$ in $Tr_{d}(\mathcal{A}_\mathcal{H})$ which is expressed by Equation (\ref{waaa2}), that is
{\small
\begin{equation}\label{A1}
\mu_{d}(j)=(m-1)^{n-1}\sum \limits_{d_{1}+\cdots+d_{n}=d}\prod \limits_{i=1}^{n}\frac{1}{(d_{i}(m-1))!}(\sum \limits_{\alpha_{i}\in [n]^{m-1}}h_{i\alpha_{i}}\frac{\partial}{\partial a_{i\alpha_{i}}})^{d_{i}}(A^{d(m-1)})_{jj}.
\end{equation}
}
Thus,
\begin{equation}\label{py}
\sum_{j=1}^n\mu_{d}(j)=Tr_{d}(\mathcal{A}_\mathcal{H}).
\end{equation}

\begin{defi}
For an $m$-uniform hypergraph $\mathcal{H}$ with $n$ vertices,
\begin{equation}\label{py11}
\sum \limits_{d=0}^{\infty}\frac{\mu_{d}(j)}{d!}
\end{equation}
is called the subgraph centrality of a vertex $j$ in $\mathcal{H}$, denoted by $C(j)$, where $\mu_{d}(j)$ is given by Equation (\ref{A1}).
\end{defi}

%We define
%\begin{equation}\label{py11}
%\sum \limits_{d=0}^{\infty}\frac{\mu_{d}(j)}{d!}
%\end{equation}
%as the  \textit{subgraph centrality} of a vertex $j$ in $\mathcal{H}$, denoted by $C(j)$.
%by Equation (\ref{py}) and Equation (\ref{py11})
By Equation (\ref{A1}), we know $C(j)$ is a real number. When $\mathcal{H}$ is $2$-uniform, since $\mu_{d}(j)=(\mathcal{A}_\mathcal{H}^{d})_{jj}$ \cite{2011Analogue}, $C(j)$ is the subgraph centrality of a graph as in Equation (\ref{wa1}).

Since the $d$th order trace of the adjacency tensor $Tr_d(\mathcal{A}_\mathcal{H})=\sum_{i=1}^k\lambda_i^d$ \cite{hu2013determinants} and Equation (\ref{py}),
\begin{align*}\label{E2}
\sum_{j=1}^n C(j)
=\sum \limits_{d=0}^{\infty}\sum_{j=1}^n\frac{\mu_{d}(j)}{d!}=\sum \limits_{d=0}^{\infty}\frac{Tr_d(\mathcal{A}_\mathcal{H})}{d!}=\sum_{i=1}^k\sum \limits_{d=0}^{\infty}\frac{\lambda_i^d}{d!}=\sum_{i=1}^k e^{\lambda_i}=EE(\mathcal{H}),
\end{align*}
where $\lambda_1,\lambda_{2},\ldots,\lambda_k$ are all the eigenvalues of $\mathcal{A}_\mathcal{H}$.

\section{Bounds for the Estrada index of hypergraphs}

The spectrum of an $m$-uniform hypergraph is said to be \textit{$m$-symmetric} if this spectrum is invariant under a rotation of an angle $2\pi/m$ in the complex plane \cite{shao2015some}. In this section, for a $3$-uniform hypergraph $\mathcal{H}$ whose spectrum is $3$-symmetric, we give an upper bound for the Estrada index in terms of energy of $\mathcal{H}$. And for $m$-uniform hypergraphs, we establish some bounds for the Estrada index.

%An $m$-uniform hypergraph with $n$ vertices is said to be complete if any $m$ vertices in $[n]$ is an edge. Let $\mathcal{K}_{n}$ be the complete $m$-uniform hypergraph with $n$ vertices and $\overline{\mathcal{K}_{n}}$ its (edgeless) complement. In this paper, $\overline{\mathcal{K}_{n}}$ is called empty hypergraph.

 The spectra of $m$-uniform power hypergraphs and hypertrees have attracted extensive attention \cite{FAN2021112329,2017On}. Their spectra are all $m$-symmetric \cite{cooper2012spectra,shao2015some,2017On,Hu2013Cored}.

\begin{lem}\label{L3}\cite{shao2015some}
Let $\mathcal{H}$ be an $m$-uniform hypergraph whose spectrum is $m$-symmetric. If $m\nmid d$, then $Tr_{d}(\mathcal{A}_\mathcal{H})=0$.
\end{lem}
Let $\mathcal{H}$ be an $m$-uniform hypergraph with $n$ vertices and $\lambda_{1},\lambda_{2},\ldots,\lambda_{k}$ be all the eigenvalues of $\mathcal{A}_\mathcal{H}$.
In this paper, $\sum_{j=1}^{k}|\lambda_{j}|$ is called the \textit{energy} of $\mathcal{H}$, denoted by $\mathcal{E}(\mathcal{H})$. When $m=2$, bounds for the Estrada index of a graph $\mathcal{H}$ were given by energy of $\mathcal{H}$ \cite{DELAPENA200770,1978The,BAMDAD2010739}.

%, where  $|\lambda_{j}|$ is modulus of $\lambda_{j}, j=1,2,\ldots, k$ $\sum_{j=1}^{n}|\lambda_{j}|$ is called \textit{energy} of a graph
For a $3$-uniform hypergraph $\mathcal{H}$ whose spectrum is $3$-symmetric, we establish an upper bound by energy of $\mathcal{H}$ for the Estrada index.
\begin{thm}\label{py3}
Let $\mathcal{H}$ be a $3$-uniform hypergraph  with $n$ vertices and $q$ edges ($q\geq 1$). If the spectrum of $\mathcal{H}$ is $3$-symmetric, then
\begin{equation*}\label{g1}
EE(\mathcal{H})\leq\frac{2(\cosh\rho-1)}{3\rho}\mathcal{E}(\mathcal{H})+k,
\end{equation*}
where $\rho$ is the spectral radius of $\mathcal{A}_{\mathcal{H}}$ and $k=2^{n-1}n$.
%The equality of the inequality is attained if and only if $\mathcal{H}$ is an empty hypergraph.
\end{thm}
\begin{proof}
%Let $\lambda_{j}=\alpha_{j}+i\beta_{j}$ be the eigenvalue of $\mathcal{A}_\mathcal{H}$, where $\alpha_{j},\beta_{j}\in \mathbb{R}$,$j\in [k]$ and $k=n(m-1)^{n-1}$. Then
%$$EE(\mathcal{H})=\sum\limits_{j=1}^{k}e^{\lambda_{j}}\leq\sum\limits_{j=1}^{k}e^{\alpha_{j}}$$
It follows from Lemma \ref{L3} that
\begin{equation}\label{o1}
EE(\mathcal{H})=\sum\limits_{d=0}^{\infty}\frac{Tr_d(\mathcal{A}_\mathcal{H})}{d!}
=\sum\limits_{l=0}^{\infty}\frac{Tr_{3l}(\mathcal{A}_\mathcal{H})}{(3l)!}
=\sum\limits_{l=0}^{\infty}\sum\limits_{j=1}^{k}\frac{\lambda_{j}^{3l}}{(3l)!}
\leq\sum\limits_{j=1}^{k}\sum\limits_{l=0}^{\infty}\frac{|\lambda_{j}|^{3l}}{(3l)!},
\end{equation}
where $\lambda_{1},\lambda_{2},\ldots,\lambda_{k}$ are all the eigenvalues of $\mathcal{A}_\mathcal{H}$.

Let $S(x)=\sum_{l=0}^{\infty}\frac{x^{3l}}{(3l)!}$, where $x$ is a real variable. We have
%$\frac{dS}{dx}=\sum_{l=1}^{\infty}\frac{x^{3l-1}}{(3l-1)!}$, $\frac{d^{2}S}{dx^{2}}=\sum_{l=1}^{\infty}\frac{x^{3l-2}}{(3l-2)!}$ and
$\frac{d^{3}S}{dx^{3}}=\sum_{l=1}^{\infty}\frac{x^{3l-3}}{(3l-3)!}
=\sum_{l=0}^{\infty}\frac{x^{3l}}{(3l)!}$.

Hence, we get the differential equation $\frac{d^{3}S}{dx^{3}}-S=0$ satisfying the conditions $S(0)=1$, $\frac{dS}{dx}\mid_{x=0}=0$ and $\frac{d^{2}S}{dx^{2}}\mid_{x=0}=0$.
Then $S(x)=\frac{2}{3}e^{-\frac{x}{2}}\cos(\frac{\sqrt{3}}{2}x)+\frac{1}{3}e^{x}$.

%Let $x=|\lambda_{j}|$.
We have
\begin{align*}
S(|\lambda_{j}|)&=\frac{2}{3}e^{-\frac{|\lambda_{j}|}{2}}\cos(\frac{\sqrt{3}}{2}|\lambda_{j}|)+\frac{1}{3}e^{|\lambda_{j}|}\\
&\leq \frac{2}{3}\frac{e^{-|\lambda_{j}|}+(\cos(\frac{\sqrt{3}}{2}|\lambda_{j}|))^{2}}{2}+\frac{1}{3}e^{|\lambda_{j}|}\\
&\leq \frac{1}{3}(e^{-|\lambda_{j}|}+e^{|\lambda_{j}|})+\frac{1}{3}=\frac{2}{3}\cosh|\lambda_{j}|+\frac{1}{3}.
\end{align*}

Let $f(x)=\frac{\cosh x-1}{x}$, $x\in(0,\rho]$. We have $\frac{d f(x)}{dx}=\frac{x\sinh x-\cosh x+1}{x^{2}}$. Let $g(x)=x\sinh x-\cosh x+1$, $x\in[0,\rho]$. We have $\frac{d g(x)}{dx}=x\cosh x\geq0$ and $\frac{d g(x)}{dx}=0$ if and only if $x=0$. Thus, $g(x)$ is strictly monotone increasing function. For $0<x\leq\rho$, we have $g(x)>g(0)=0$. Hence, $\frac{d f(x)}{dx}>0$, $x\in(0,\rho]$. Then $f(x)$ is  strictly monotone increasing function. So $\frac{\cosh x-1}{x}\leq\frac{\cosh\rho-1}{\rho}$, $x\in(0,\rho]$, that is $\cosh x\leq\frac{\cosh\rho-1}{\rho}x+1$, $x\in(0,\rho]$. When $x=0$, the above inequality obviously holds. Thus, $\cosh x\leq\frac{\cosh\rho-1}{\rho}x+1$, $x\in[0,\rho]$.

%When $x=0$, we have $\cosh x=\frac{\cosh\rho-1}{\rho}x+1$.

%Since $\cosh x\leq\frac{\cosh\rho-1}{\rho}x+1, x\in[0,\rho]$, we have
So,
\begin{align}\label{pf1}
S(|\lambda_{j}|)&\leq\frac{2}{3}\cosh|\lambda_{j}|+\frac{1}{3}\notag\\
&\leq \frac{2}{3}(\frac{\cosh\rho-1}{\rho}|\lambda_{j}|+1)+\frac{1}{3}=\frac{2(\cosh\rho-1)}{3\rho}|\lambda_{j}|+1.
\end{align}
By Equation (\ref{o1}) and (\ref{pf1}), we have
\begin{align*}
EE(\mathcal{H})&\leq\sum_{j=1}^{k}S(|\lambda_{j}|)\\
&\leq\frac{2(\cosh\rho-1)}{3\rho}\sum_{j=1}^{k}|\lambda_{j}|+k=\frac{2(\cosh\rho-1)}{3\rho}\mathcal{E}(\mathcal{H})+k.
\end{align*}

\end{proof}

%Similar the proof of Theorem \ref{py3}, we easily obtain the following result.
%\begin{cor}
%Let $\mathcal{H}$ be a $4$-uniform hypergraph  with $n$ vertices. If the spectrum of $\mathcal{A}_\mathcal{H}$ is $4$-symmetric, then
%\begin{equation}\label{g1}
%EE(\mathcal{H})\leq\frac{\cosh\rho-1}{2\rho}\sum_{j=1}^{k}|\lambda_{j}|+k,
%\end{equation}
%where $\rho$ is the spectral radius of $\mathcal{A}_{\mathcal{H}}$ and $\lambda_{1},\lambda_{2},...,\lambda_{k}$ are all the eigenvalues of $\mathcal{A}_\mathcal{H}$, $k=3^{n-1}n$.
%%The equality of the inequality is attained if and only if $\mathcal{H}$ is an empty hypergraph.
%\end{cor}

The following is a lower bound for the Estrada index of $m$-uniform hypergraphs.
\begin{thm}\label{wa3}
Let $\mathcal{H}$ be an $m$-uniform hypergraph with $n$ vertices and $q$ edges. Then
\begin{equation}\label{0}
EE(\mathcal{H})\geq\frac{qm^{m-2}(m-1)^{n-m-1}}{(m-2)!}+n(m-1)^{n-1},
\end{equation}
%where $\rho$ is the spectral radius of $\mathcal{A}_\mathcal{H}$.
%The equalities on both sides of the above inequality are attained respectively if and only if $\mathcal{H}$ is an empty hypergraph.
equality holds if and only if $\mathcal{H}$ is an empty hypergraph.
\end{thm}

\begin{proof}
From Lemma \ref{L2}, we have
%$Tr_{0}(\mathcal{A}_\mathcal{H})=n(m-1)^{n-1}$ and $Tr_{d}(\mathcal{A}_\mathcal{H})=0$, $d=1,2,\ldots,m-1$. Then
\begin{align*}
   EE(\mathcal{H})&=\sum \limits_{d=0}^{\infty}\frac{Tr_{d}(\mathcal{A}_\mathcal{H})}{d!}=n(m-1)^{n-1}+\sum \limits_{d=m}^{\infty}\frac{Tr_{d}(\mathcal{A}_\mathcal{H})}{d!}\\ %=n(m-1)^{n-1}+\sum \limits_{d=m}^{\infty}\frac{Tr_{d}(\mathcal{A}_\mathcal{H})}{d!} \\
  &\geq n(m-1)^{n-1}+\frac{Tr_{m}(\mathcal{A}_\mathcal{H})}{m!}=n(m-1)^{n-1}+\frac{qm^{m-2}(m-1)^{n-m-1}}{(m-2)!}.
\end{align*}
%It follows from Lemma \ref{L2} (3) that
%$$
%EE(\mathcal{H})\geq n(m-1)^{n-1}+\frac{qm^{m-2}(m-1)^{n-m-1}}{(m-2)!}.
%$$

When $\mathcal{H}$ is an empty hypergraph, all eigenvalues of $\mathcal{A}_\mathcal{H}$ are 0. It is easy to see $EE(\mathcal{H})=\sum_{i=1}^{k}e^{0}=k=n(m-1)^{n-1}$.
%Then the equalities on both sides of Inequality (\ref{0}) are attained, respectively.
Then the equality of Inequality (\ref{0}) holds.

On the other hand, if the equality of Inequality (\ref{0}) holds, it follows from the proof of the above inequality that $Tr_{d}(\mathcal{A}_\mathcal{H})=0$ for all $d\geq m+1$, that is
\begin{equation}\label{wp}
\sum_{j=1}^{k}\lambda_{j}^{d}=0,~  \text{for all} ~d\geq m+1,
\end{equation}
where $\lambda_{1},\lambda_{2},\ldots,\lambda_{k}$ are all the eigenvalues of $\mathcal{A}_\mathcal{H}$.
Without loss of generality,
suppose that $\lambda_1,\lambda_2,\dots, \lambda_s$ are all the distinct eigenvalues among $\lambda_1,\lambda_2,\dots, \lambda_k$, and $l_i\geq1$ is the
multiplicity of $\lambda_i$, $i=1,2,\ldots,s$.

If $s=1$, then all the eigenvalues of $\mathcal{A}_\mathcal{H}$ are the same. Since $Tr_1(\mathcal{A}_\mathcal{H})=0$, all the eigenvalues of $\mathcal{A}_\mathcal{H}$ are $0$.

If $s\geq2$, let
$M=\left(
\begin{array}{cccc}
1 & 1 & \cdots & 1 \\
\lambda_1 & \lambda_2 & \cdots  & \lambda_s \\
\cdots& \cdots & \cdots & \cdots \\
\lambda_1^{s-1} & \lambda_2^{s-1} & \cdots  & \lambda_s^{s-1} \\
\end{array}
\right)
$. By Equation (\ref{wp}), we have $M(l_1\lambda_1^{m+1},l_2\lambda_2^{m+1},\dots,l_s\lambda_s^{m+1})^{\rm T}=0$. Since $\det(M)\neq 0$, $(l_1\lambda_1^{m+1},l_2\lambda_2^{m+1},\dots,\\l_s\lambda_s^{m+1})^{\rm T}=0$, that is $\lambda_{j}=0$, $j=1,2,\ldots,s$. It contradicts $s\geq2$.

So all the eigenvalues of $\mathcal{A}_\mathcal{H}$ are $0$. Then the spectral radius of $\mathcal{A}_\mathcal{H}$ is $0$. The spectral radius is greater than or equal to the average degree of $\mathcal{H}$ \cite{cooper2012spectra}. Hence, the average degree is equal to 0,
%Without loss of generality, suppose that $\mathcal{H}$ is connected. Perron-Frobenius Theorem for tensors\cite{21} gives there exists an  $x>0$ such that $\mathcal{A}_\mathcal{H} x^{m-1}=\rho x^{[m-1]}=0$. Hence, $\mathcal{A}_\mcooper2012spectra,athcal{H}$ is a zero tensor, that is
that is $\mathcal{H}$ is an empty hypergraph.

%If the equality on the right hand side of Inequality (\ref{0}) holds, then all eigenvalues of $\mathcal{H}$ are equal to the spectral radius $\rho$. And since $Tr_{1}(\mathcal{A}_\mathcal{H})=0$, we have $\rho=0$. Hence, $\mathcal{H}$ is an empty hypergraph.

Thus, equality holds in Inequality (\ref{0}) if and only if $\mathcal{H}$ is an empty hypergraph.
\end{proof}

Next we establish an upper bound for the Estrada index of $m$-uniform hypergraphs with $m\geq 3$.
\begin{thm}\label{3}
Let $\mathcal{H}$ be an $m$-uniform hypergraph with $n$ vertices, $q$ edges and $m\geq 3$. Then

%(1)
%\begin{equation}\label{11}
% EE(\mathcal{H})\leq k-1+e^{\sqrt{r}};
%\end{equation}
%
%(2)
\begin{equation}\label{11}
 EE(\mathcal{H})\leq k-1+e^{r}+\frac{qm^{m-2}(m-1)^{n-m-1}}{(m-2)!}-\sum_{d=1}^{m}\frac{r^{d}}{d!},
\end{equation}
where $r=\sqrt{2\sum_{j=1}^{k}(Re(\lambda_{j}))^{2}}$, $\lambda_{1}, \lambda_{2},\ldots,\lambda_{k}$ are all the eigenvalues of $\mathcal{A}_\mathcal{H}$, $k=n(m-1)^{n-1}$. Equality holds if and only if $\mathcal{H}$ is an empty hypergraph.
\end{thm}
\begin{proof}
Let $\lambda_{j}=\alpha_j+\textbf{i}\beta_j$, $\alpha_j, \beta_j\in \mathbb{R}$, $j=1,2,\ldots,k$, $\textbf{i}^{2}=-1$. From $Tr_{2}(\mathcal{A}_\mathcal{H})=\sum_{j=1}^{k}\lambda_{j}^{2}$ and Lemma \ref{L2} (2), we have
$$
Tr_{2}(\mathcal{A}_\mathcal{H})
=\sum\limits_{j=1}^{k}(\alpha_{j}^{2}-\beta_{j}^{2})+\textbf{i}\sum\limits_{j=1}^{k}2\alpha_{j}\beta_{j}=0.
$$
%Note that $Tr_{2}(\mathcal{A}_\mathcal{H})$ is a nonnegative real number,
Then
$$
\sum\limits_{j=1}^{k}(\alpha_{j}^{2}-\beta_{j}^{2})=0.$$
So
%$$
%\sum\limits_{j=1}^{k}\beta_{j}^{2}=\sum\limits_{j=1}^{k}\alpha_{j}^{2}-Tr_{2}(\mathcal{A}_\mathcal{H})
%$$
%which yields
\begin{align}\label{13}
\sum \limits_{j=1}^{k}\left| \lambda_{j}\right|^{2}=\sum \limits_{j=1}^{k}(\alpha_{j}^{2}+\beta_{j}^{2})=2\sum\limits_{j=1}^{k}\alpha_{j}^{2}.
\end{align}

Note that  $EE(\mathcal{H})$, $Tr_{d}(\mathcal{A}_\mathcal{H})$ are nonnegative real numbers, $d=0,1,2, \ldots,$ and $\sum_{d=0}^{^{\infty}}\frac{\left| \lambda_{j}\right|^{d}}{d!}$ is convergent.
%$\sum_{i=1}^{n}a_{i}^{k}\leq(\sum_{i=1}^{n}a_{i}^{2})^{^{k/2}}$, where $a_{i}\geq0,i\in[n]$, $k\geq 2$.
We have
\begin{align*}
EE(\mathcal{H})&=\sum\limits_{d=0}^{m}\frac{Tr_{d}(\mathcal{A}_\mathcal{H})}{d!}+\sum\limits_{d=m+1}^{\infty}\frac{\sum_{j=1}^{k}\lambda_{j}^{d}}{d!}\\
%&=\sum\limits_{d=0}^{m}\frac{Tr_{d}(\mathcal{A}_\mathcal{H})}{l!}+|\sum\limits_{l=m+1}^{\infty}\frac{\sum_{j=1}^{k}\lambda_{j}^{l}}{l!}|\\
&\leq\sum_{d=0}^{m}\frac{Tr_{d}(\mathcal{A}_\mathcal{H})}{d!}+\sum \limits_{d=m+1}^{\infty}\frac{1}{d!}\sum \limits_{j=1}^{k}\left| \lambda_{j}\right|^{d}.\\
%=\sum_{l=0}^{m}\frac{Tr_{l}(\mathcal{A}_\mathcal{H})}{l!}+\sum \limits_{l=m+1}^{\infty}\frac{1}{l!}\sum \limits_{j=1}^{k}(\left| \lambda_{j}\right|^{2})^{\frac{l}{2}}\\
%&\leq\sum_{l=0}^{m}\frac{Tr_{l}(\mathcal{A}_\mathcal{H})}{l!}+\sum \limits_{l=m+1}^{\infty}\frac{1}{l!}(\sum \limits_{j=1}^{k}\left| \lambda_{j}\right|^{2})^{\frac{l}{2}}.
\end{align*}
%=&\sum_{l=0}^{m}\frac{Tr_{l}(\mathcal{A}_\mathcal{H})}{l!}+\sum \limits_{j=1}^{k}\sum \limits_{l=m+1}^{\infty}\frac{(\cos(\theta_{j} l)+\textbf{i}\sin(\theta_{j} l))\left| \lambda_{j}\right|^{l}}{l!}.
%\end{align*}
%Note that $EE(\mathcal{H})$ is a nonnegative real number. Then
%\begin{align*}
%EThe equality in the inequality is attained if and only if $\mathcal{H}$ is an empty hypergraph.

                                               %E(\mathcal{H})&=\sum_{l=0}^{m}\frac{Tr_{l}(\mathcal{A}_\mathcal{H})}{l!}+\sum \limits_{j=1}^{k}\sum \limits_{l=m+1}^{\infty}\frac{\cos(\theta_{j} l)\left| \lambda_{j}\right|^{l}}{l!}\\
%&\leq\sum_{l=0}^{m}\frac{Tr_{l}(\mathcal{A}_\mathcal{H})}{l!}+\sum \limits_{j=1}^{k}\sum \limits_{l=m+1}^{\infty}\frac{\left|\cos(\theta_{j} l)\right|\left| \lambda_{j}\right|^{l}}{l!}\\
%\end{align*}
For $\left| \lambda_{j}\right|, j\in[k]$ and integer $d\geq2$, by Cauchy-Schwarz Inequality, we have
\begin{align*}
(\sum \limits_{j=1}^{k}\left| \lambda_{j}\right|^{d})^{2}&=(\sum \limits_{j=1}^{k}(\left| \lambda_{j}\right|^{d-1}\left| \lambda_{j}\right|))^{2}\\
&\leq\sum \limits_{j=1}^{k}\left| \lambda_{j}\right|^{2(d-1)}\sum \limits_{j=1}^{k}\left| \lambda_{j}\right|^{2}
=\sum \limits_{j=1}^{k}(\left| \lambda_{j}\right|^{2})^{(d-1)}\sum \limits_{j=1}^{k}\left| \lambda_{j}\right|^{2}\\
&\leq (\sum \limits_{j=1}^{k}\left| \lambda_{j}\right|^{2})^{d-1}\sum \limits_{j=1}^{k}\left| \lambda_{j}\right|^{2}=(\sum \limits_{j=1}^{k}\left| \lambda_{j}\right|^{2})^{d},
\end{align*}
that is $\sum_{j=1}^{k}\left| \lambda_{j}\right|^{d}\leq(\sum_{j=1}^{k}\left| \lambda_{j}\right|^{2})^{\frac{d}{2}}$. So,
$$
EE(\mathcal{H})\leq\sum_{d=0}^{m}\frac{Tr_{d}(\mathcal{A}_\mathcal{H})}{d!}+\sum \limits_{d=m+1}^{\infty}\frac{1}{d!}(\sum \limits_{j=1}^{k}\left| \lambda_{j}\right|^{2})^{\frac{d}{2}}.
$$
It follows from Equation (\ref{13}) that
\begin{align*}
EE(\mathcal{H})&\leq \sum_{d=0}^{m}\frac{Tr_{d}(\mathcal{A}_\mathcal{H})}{d!}+\sum \limits_{d=m+1}^{\infty}\frac{1}{d!}(2\sum\limits_{j=1}^{k}\alpha_{j}^{2})^{\frac{d}{2}}\\
  &=\sum_{d=0}^{m}\frac{Tr_{d}(\mathcal{A}_\mathcal{H})}{d!}+e^{\sqrt{2\sum\limits_{j=1}^{k}\alpha_{j}^{2}}}-\sum_{d=0}^{m}\frac{(2\sum\limits_{j=1}^{k}\alpha_{j}^{2})^{\frac{d}{2}}}{d!}.
\end{align*}
By Lemma \ref{L2}, we have
$$
EE(\mathcal{H})\leq k-1+e^{\sqrt{2\sum\limits_{j=1}^{k}\alpha_{j}^{2}}}+\frac{Tr_{m}(\mathcal{A}_\mathcal{H})}{m!}-\sum_{d=1}^{m}\frac{(2\sum\limits_{j=1}^{k}\alpha_{j}^{2})^{\frac{d}{2}}}{d!}
,$$
where $Tr_{m}(\mathcal{A}_\mathcal{H})=qm^{m-1}(m-1)^{n-m}$. Thus, we obtain Inequality (\ref{11}).

If $\mathcal{H}$ is an empty hypergraph,  all eigenvalues of $\mathcal{A}_\mathcal{H}$ are 0. It is easy to see $EE(\mathcal{H})=\sum_{i=1}^{k}e^{0}=k$. Then the equality in Inequality (\ref{11}) holds.

On the other hand, suppose the equality in Inequality (\ref{11}) holds. From the proof of the above inequality, we have $\sum_{d=m+1}^{\infty}\frac{1}{d!}(\sum_{j=1}^{k}|\lambda_{j}|^{d}-\sum_{j=1}^{k}\lambda_{j}^{d})=0$ and note that $\sum_{j=1}^{k}\lambda_{j}^{d}$ is a nonnegative real number, $d=m+1,m+2, \ldots.$ Then
\begin{equation*}
\sum_{j=1}^{k}\lambda_{j}^{d}=\sum_{j=1}^{k}|\lambda_{j}|^{d},d=m+1,m+2, \ldots.
\end{equation*}
Let $\lambda_{j}=|\lambda_{j}|e^{\textbf{i}\theta_{j}},j=1,2,\ldots,k$. We have $$\sum\limits_{j=1}^{k}|\lambda_{j}|^{d}(1-\cos(\theta_{j}d))=0, d=m+1,m+2, \ldots.$$
Then
$|\lambda_{j}|^{d}(1-\cos(\theta_{j}d))=0$, $d=m+1,m+2, \ldots, \text{and}~ j=1,2,\ldots,k.$
Thus, we have $\lambda_{j}=0$ or $\cos(\theta_{j}d)=1$, $d=m+1,m+2, \ldots,$ for each $j\in[k]$. If there exists $j\in[k]$ such that $\cos(\theta_{j}d)=1$, $d=m+1,m+2, \ldots.$ Then we have $\theta_{j}d=2l_{d}\pi,$ where $l_{d}$ is an integer and $d=m+1,m+2, \ldots.$ So,
$$
\theta_{j}=\frac{l_{m+1}}{m+1}2\pi=\frac{l_{m+2}}{m+2}2\pi=\cdots.
$$
Let $
\frac{l_{m+1}}{m+1}=\frac{l_{m+2}}{m+2}=\cdots=t.$ Then
$$
l_{m+2}=t(m+2)=t(m+1)+t=l_{m+1}+t.
$$
Since $l_{m+1}, l_{m+2}$ are integers, $t$ is an integer.  So, $\theta_{j}$ is equal to an integer multiple of $2\pi$. Hence, $\lambda_{j}$ is a nonnegative  real number.
%we get $\theta_{j}=0$ and $l_{d}=0, d=m+1,m+2, \ldots, $ that is $\lambda_{j}$ is a nonnegative  real number.
Since $Tr_1(\mathcal{A}_\mathcal{H})=0$,
all eigenvalues of $\mathcal{A}_\mathcal{H}$ are $0$.
In the proof of Theorem \ref{wa3}, we proved that $\mathcal{H}$ is an empty hypergraph when all eigenvalues of $\mathcal{A}_\mathcal{H}$ are $0$. Thus,
we get $\mathcal{H}$ is an empty hypergraph. Therefore, the equality holds in Inequality (\ref{11}) if and only if $\mathcal{H}$ is an empty hypergraph. %Similarly, the equality in Inequality (\ref{12}) is attained if and only if $\mathcal{H}$ is an empty hypergraph.
%If there exists $j\in[k]$, such that $\cos(\theta_{j}d)=1$. Then we have $\theta_{j}d=2l\pi,$ where $d=m+1,m+2, \ldots,$ and $l$ is an integer. So, we get $\theta_{j}=0$ and $l=0$, that is $\lambda_{j}=|\lambda_{j}|$.
\end{proof}
\section{Expressions for the subgraph centrality of hypergraphs}

In this section, for a uniform hypergraph $\mathcal{H}$, we give two expressions of the subgraph centrality by the number of  Eulerian closed walks of the multi-digraphs  associated with $\mathcal{H}$ and the number of arborescences of the multi-digraphs associated with  $\mathcal{H}$, respectively. The explicit expression of $\mu_{d}(j)$ of $\mathcal{H}$ is given for $d=0,1,2,\ldots,m$.
%
%For an $m$-uniform hypergraph $\mathcal{H}$, we give an expression for $\mu_d(j)$ by using the number of the Eulerian closed walks of the multi-digraphs  associated with  $\mathcal{H}$. When $m=2$, $\mu_{d}(j)=(\mathcal{A}_\mathcal{H}^{d})_{jj}$ is equal to the number of closed walks of length $d$ starting and ending at the vertex $j$ in $\mathcal{H}$ \cite{30,1+}. Thus, when $m=2$, $\mu_d(j)$ which is attained by the number of the Eulerian closed walks of the multi-digraphs  associated with  $\mathcal{H}$ is equal to the number of closed walks of length $d$ starting and ending at the vertex $j$ in $\mathcal{H}$.

%For a positive integer $n$,

%where $[n]=\{1,2,\cdots,n\}$.
%For $(i_1,i_2,\cdots,i_m)\in[n]^{m}$, it is also writed by $i_1i_2\cdots i_m$.
For an integer $d>0$, let $\alpha=i_1\cdots i_{d}\in[n]^{d}$.
If $\alpha$ satisfies $i_1\leq \cdots \leq i_{d}$, then $\alpha$ is called \textit{ascending order}.
Let
$$
\mathcal{F}_{d}=\{(i_{1}\alpha_{1},\ldots,i_{d}\alpha_{d})|~1\leq i_{1}\leq\cdots\leq i_{d}\leq n,~ \alpha_{k}\in[n]^{m-1}, ~k=1,\ldots,d\}.
$$

%For two vertices  $v_{1}, v_{2}$ of a digraph, let $(v_{1},v_{2})$ denote the arc from vertex $v_{1}$ to $v_{2}$.

Let $F=(i_{1}\alpha_{1},\ldots,i_{d}\alpha_{d})\in \mathcal{F}_{d}$, where $i_{k}\alpha_{k}=i_{k}j_1^{(k)}\cdots j_{m-1}^{(k)}$, $k=1,2,\ldots,d$. Let the set of arcs from $i_{k}$ to $j_1^{(k)}, j_2^{(k)}, \ldots, j_{m-1}^{(k)}$ be $$E(i_{k}\alpha_{k})=\{(i_{k},j_1^{(k)}),(i_{k},j_2^{(k)}),\ldots,(i_{k},j_{m-1}^{(k)})\}. $$ Let the arc multi-set $$\widetilde{E}(F)=\bigcup_{k=1}^{d}E(i_k\alpha_k).$$

Let multi-digraph $D(F)=(V(\widetilde{E}(F)), \widetilde{E}(F)).$ An \textit{Eulerian tour} of  the multi-digraph $D(F)$ is a closed walk that traverses each arc of $D(F)$ exactly once \cite{matthew2016multi}. Let $\epsilon(F,j)$ be the set of all Eulerian tours starting and ending at the vertex $j$ in $D(F)$.
For any arc $(i,j)$ from $i$ to $j$ in $D(F)$, let $\omega(i,j)$ be multiplicity of arc $(i,j)$.
Let $\widehat{D}(F)$ be the digraph formed by removing duplicate arcs of $D(F)$.
%arc from $i$ to $j$ in $D(F)$.
In this paper, an \textit{Eulerian closed walk} of the multi-digraph $D(F)$ is a closed walk that uses each arc $(i,j)$ of $\widehat{D}(F)$  exactly $\omega(i,j)$ times. Let $\textbf{W}(F,j)$ be the set of all Eulerian closed walks of $D(F)$ starting and ending at the vertex $j$.
The multi-digraph $D(F)$ is called \textit{Eulerian} if $D(F)$ has Eulerian closed walks.
%Note that duplicate arcs of $D(F)$ are distinguished. If duplicate arcs of $D(F)$ are not distinguished, sequence of above alternative vertices and arcs is called the \textit{Eulerian closed walk}.

%\begin{example}
%\begin{figure}[h]
%\centering
%\includegraphics[scale=0.4]{3}
%\caption{The multi-digraph $D$}
%\end{figure}
%Let $D$ be the multi-digraph shown in Figure 1. There are 8 Eulerian tours, namely $v_{1}e_{1}v_{2}e_{3}v_{1}e_{2}v_{2}e_{4}v_{1},$ $v_{1}e_{1}v_{2}e_{4}v_{1}e_{2}v_{2}e_{3}v_{1}, $ $v_{1}e_{2}v_{2}e_{3}v_{1}e_{1}v_{2}e_{4}v_{1},$  $v_{1}e_{2}v_{2}e_{4}v_{1}e_{1}v_{2}e_{3}v_{1},$  $v_{2}e_{3}v_{1}e_{2}v_{2}e_{4}v_{1}e_{1}v_{2},$  $v_{2}e_{3}v_{1}e_{1}v_{2}e_{4}v_{1}e_{2}v_{2},$ $v_{2}e_{4}v_{1}e_{2}v_{2}e_{3}v_{1}e_{1}v_{2},$ $v_{2}e_{4}v_{1}e_{1}v_{2}e_{3}v_{1}e_{2}v_{2}.$ There are 2 Eulerian closed walks, namely $v_{1}(v_{1},v_{2})v_{2}(v_{2},v_{1})\\v_{1}(v_{1},v_{2})v_{2}(v_{2},v_{1})v_{1},$ $v_{2}(v_{2},v_{1})v_{1}(v_{1},v_{2})v_{2}(v_{2},v_{1})v_{1}(v_{1},v_{2})v_{2}$.
%\end{example}

For $F=( i_{1}\alpha_{1},\ldots,i_{d}\alpha_{d})\in \mathcal{F}_{d}$,
%$a_{i_{k}\alpha_{k}} $ is an auxiliary variable ($k=1,\cdots,d$), and
let the differential operator $\partial(F):=\prod_{k=1}^{d}\frac{\partial}{\partial a_{i_{k}\alpha_{k}}}$.
%where $\frac{\partial}{\partial a_{i_{k}\alpha_{k}}}:=\prod\limits_{s=1}^{m-1}\frac{\partial}{\partial a_{i_{k}j^{(k)}_{s}}}$, $\alpha_{k}=j_{1}^{(k)}\cdots j_{m-1}^{(k)}$.

\begin{lem}\label{z1}
Let $F\in\mathcal{F}_{d}$. Then
\begin{equation*}
\partial(F)(A^{d(m-1)})_{jj}=b(F)|\textbf{W}(F,j)|,
\end{equation*}
where $A $ is an auxiliary matrix of order $n$, and $b(F)$ is the product of the factorials of the multiplicities of all the arcs in $D(F)$.
\end{lem}

\begin{proof}
For the matrix $A = (a_{ij})$, we have
%let $i_{k}\alpha_{k}=i_{k}j_{1}^{(k)}\cdots j_{m-1}^{(k)}\in[n]^m, k=1,2,...,d$. For
%$$
%(A^{d(m-1)})_{jj}=\sum\limits_{i_{2},\cdot\cdot\cdot,i_{d(m-1)}=1}^{n}a_{ji_{2}}a_{i_{2}i_{3}}\cdot\cdot\cdot a_{i_{d(m-1)}j}
%$$
%and
\begin{align*}
\partial(F)(A^{d(m-1)})_{jj}&=\sum\limits_{l_{2},\ldots,l_{d(m-1)}=1}^{n}\partial(F)(a_{jl_{2}}a_{l_{2}l_{3}}\cdots a_{l_{d(m-1)}j}).
%&=\sum\limits_{l_{2},\ldots,l_{d(m-1)}=1}^{n}\prod_{k=1}^{d}\prod\limits_{s=1}^{m-1}\frac{\partial}{\partial a_{i_{k}j^{(k)}_{s}}}(a_{jl_{2}}a_{l_{2}l_{3}}\cdot\cdot\cdot a_{l_{d(m-1)}j}).
\end{align*}
%From the proof of Lemma 2.5 in \cite{17},
Since $a_{jl_{2}}a_{l_{2}l_{3}}\cdot\cdot\cdot a_{l_{d(m-1)}j}$ is a polynomial of degree $d(m-1)$, $\partial(F)(a_{jl_{2}}a_{l_{2}l_{3}}\cdot\cdot\cdot a_{l_{d(m-1)}j})\neq 0$ if and only if $\partial(F)(a_{jl_{2}}a_{l_{2}l_{3}}\cdot\cdot\cdot a_{l_{d(m-1)}j})=b(F).$
And in this case, the arc multi-set  $\{(j,l_2),(l_2,l_3),\ldots, (l_{d(m-1)},j)\}=\widetilde{E}(F)$,
that is there exists an Eulerian closed walk of $D(F)$ starting and ending at the vertex $j$.
Hence,
$$
\partial(F)(A^{d(m-1)})_{jj}=\sum\limits_{l_{2},\ldots,l_{d(m-1)}=1}^{n}\partial(F)(a_{jl_{2}}a_{l_{2}l_{3}}\cdot\cdot\cdot a_{l_{d(m-1)}j})=b(F)|\textbf{W}(F,j)|.
$$
\end{proof}

%\{(i_{1},j_{1}^{(1)}),(i_{1},j_{2}^{(1)}),\ldots,\\(i_{1},j_{m-1}^{(1)}),(i_{2},j_{1}^{(2)}),(i_{2},j_{2}^{(2)}),\ldots,(i_{2},j_{m-1}^{(2)}),\ldots,(i_{d},j_{1}^{(d)}),(i_{d},j_{2}^{(d)}),\ldots,(i_{d},j_{m-1}^{(d)})\}=

Let $\mathcal{H}=(V(\mathcal{H}),E(\mathcal{H}))$ be an $m$-uniform hypergraph with $n$ vertices, where $V(\mathcal{H})=\{1,\ldots,n\}$.
Let $\mathcal{F}_{d}(\mathcal{H})=\{(i_{1}\alpha_{1},\ldots,i_{d}\alpha_{d})\in\mathcal{F}_{d}|i_{k}\alpha_{k}\in E(\mathcal{H}), \alpha_{k}~\text{is}\\ \text{ascending order},~k=1,\ldots,d\}.$
%Clearly, $\mathcal{F}_{d}(\mathcal{H})\subset \mathcal{F}_{d}$.

Let $\mathcal{F}_{d}^{(j)}(\mathcal{H})=\{F\in\mathcal{F}_{d}(\mathcal{H})| \text{$D(F)$ is Eulerian and contains the vertex $j$}\}$.

\begin{thm}\label{z2}
Let $\mathcal{H}=(V(\mathcal{H}), E(\mathcal{H}))$ be an $m$-uniform hypergraph with $n$ vertices. Then
\begin{equation*}\label{py1}
C(j)=\sum \limits_{d=0}^{\infty}\frac{\mu_{d}(j)}{d!},
\end{equation*}
where
\begin{equation*}
\mu_{d}(j)=(m-1)^{n-1}\sum\limits_{F\in\mathcal{F}_{d}^{(j)}(\mathcal{H})}\frac{b(F)}{c(F)}|\textbf{W}(F,j)|,
\end{equation*}
$j\in V(\mathcal{H})$, $b(F)$ is the product of the factorials of the multiplicities of all the arcs in $D(F)$, $c(F)$ is the product of the factorials of the outdegrees of all the vertices in $D(F)$ and $\textbf{W}(F,j)$ is the set of all Eulerian closed walks of $D(F)$ starting and ending at the vertex $j$.
%$E(F)$ is the number of $F'\in F_{d}$ satisfying $E(F')=E(F)$.
\end{thm}

\begin{proof}
Let  $V(\mathcal{H})=\{1,\ldots,n\}$, the adjacency tensor $\mathcal{A}_\mathcal{H}=(h_{\alpha})$, $F=(i_{1}\alpha_{1},\ldots,\\i_{d}\alpha_{d})\in \mathcal{F}_{d}$
and $ \pi_{F}(\mathcal{A}_\mathcal{H})=\prod_{k=1}^{d}h_{i_{k}\alpha_{k}}$.
Using Formula (2.9) in \cite{shao2015some}
%for given $d_1,\cdots,d_n$ in Equation (\ref{A1}), we have
%$$
%\prod \limits_{i=1}^{n}\frac{1}{(d_{i}(m-1))!}(\sum \limits_{\alpha_{i}\in [n]^{m-1}}h_{i\alpha_{i}}\frac{\partial}{\partial a_{i\alpha_{i}}})^{d_{i}}=\frac{1}{c(F)}\pi_{F}(\mathcal{A}_\mathcal{H})\partial(F),
%$$
%where $F=((1,\alpha_1)^{d_1},(2,\alpha_2)^{d_2},(n,\alpha_n)^{d_n})\in \mathcal{F}_d$ and $(i,\alpha_i)^{d_i}$ is a  multi-set of $(i,\alpha_i)$ with multiplicity $d_i$, $i=1,\cdots n$. Then
$$
\sum \limits_{d_{1}+\cdots+d_{n}=d}\prod \limits_{i=1}^{n}\frac{1}{(d_{i}(m-1))!}(\sum \limits_{\alpha_{i}\in [n]^{m-1}}h_{i\alpha_{i}}\frac{\partial}{\partial a_{i\alpha_{i}}})^{d_{i}}=\sum\limits_{F\in\mathcal{F}_{d}}\frac{1}{c(F)}\pi_{F}(\mathcal{A}_\mathcal{H})\partial(F),
$$
%where $c(F)$ is the product of the factorials for the outdegrees of all the vertices in $E(F)$.
we get
$$
\mu_{d}(j)=(m-1)^{n-1}\sum\limits_{F\in\mathcal{F}_{d}}\frac{1}{c(F)}\pi_{F}(\mathcal{A}_\mathcal{H})\partial(F)(A^{d(m-1)})_{jj}.
$$
By Lemma \ref{z1}, we have
\begin{equation}\label{m1}
\mu_{d}(j)=(m-1)^{n-1}\sum\limits_{F\in\mathcal{F}_{d}}\frac{b(F)}{c(F)}\pi_{F}(\mathcal{A}_\mathcal{H})|\textbf{W}(F,j)|.
\end{equation}

According to the definition of adjacency tensors, we have
$$
\pi_{F}(\mathcal{A}_{\mathcal{H}})=\prod_{k=1}^{d}h_{i_{k}\alpha_{k}}=\begin{cases}
\frac{1}{((m-1)!)^{d}},& \text{if } ~i_{k}\alpha_{k}\in E(\mathcal{H}), k=1,\ldots,d,\notag \\
0,& \text{otherwise}.
\end{cases}
$$
Then the $F$ in Equation (\ref{m1}) such that $\pi_{F}(\mathcal{A}_{\mathcal{H}})\neq 0$ if and only if $i_{k}\alpha_{k}\in E(\mathcal{H})$, $k=1,\ldots,d$.
For $F=(i_{1}\alpha_{1},\ldots,i_{d}\alpha_{d})\in\mathcal{F}_d(\mathcal{H})$, $\alpha_{k}$ is ascending order, $k=1,\ldots,d$. Thus, each $F\in\mathcal{F}_d(\mathcal{H})$ corresponds to $((m-1)!)^d$ elements in $\mathcal{F}_d$.

%For $F=(i_{1}j_{1}^{(1)}\cdots j_{m-1}^{(1)},\ldots,i_{d}j_{1}^{(d)}\cdots j_{m-1}^{(d)})\in\mathcal{F}_d(\mathcal{H})$, let $$\mathcal{F}_{d}^\prime=\{(i_{1}\alpha_{1},\ldots,i_{d}\alpha_{d})|\text{$\alpha_k=j_{\sigma(1)}^{(k)}\cdots j_{\sigma(m-1)}^{(k)}$ for all $\sigma\in S_{m-1}$, $k=1,\ldots,d$}\}\subset \mathcal{F}_{d}.$$ Then the number of elements of $\mathcal{F}_{d}^\prime$ is $((m-1)!)^d$.
%$$\mathcal{F}_{d}^\prime=\{(i_{11}\alpha_{1},\cdot\cdot\cdot,i_{dd}\alpha_{d})|\text{$\alpha_k$ is any permutation of $\{i_{k2},\cdots,i_{km}\}$ , $k=1,\cdots,d$}\}.$$  We have $\mathcal{F}_{d}^\prime\subset \mathcal{F}_{d}$ and $|\mathcal{F}_{d}^\prime|=((m-1)!)^d.$
%there are $(m-1)!$ different ways to arrange the index array $i_{kk}\alpha_{k}$.  $\mathcal{F}_d$ satisfying $\alpha_k$ having the same elements as $\{i_{k2},\cdots,i_{km}\}$, $k=1,\cdots,d$ is $((m-1)!)^d$.
Hence,
\begin{align}\label{wei1}
\mu_{d}(j)&=(m-1)^{n-1}\sum\limits_{F\in\mathcal{F}_{d}(\mathcal{H})}((m-1)!)^d\frac{b(F)}{c(F)}\frac{1}{((m-1)!)^{d}}|\textbf{W}(F,j)|\notag\\\notag
&=(m-1)^{n-1}\sum\limits_{F\in\mathcal{F}_{d}(\mathcal{H})}\frac{b(F)}{c(F)}|\textbf{W}(F,j)|\\
&=(m-1)^{n-1}\sum\limits_{F\in\mathcal{F}_{d}^{(j)}(\mathcal{H})}\frac{b(F)}{c(F)}|\textbf{W}(F,j)|.
\end{align}
Substituting Equation (\ref{wei1}) into Equation (\ref{py11}), we obtain the expression of $C(j)$ by the number of Eulerian closed walks of the multi-digraphs associated with $\mathcal{H}$.
\end{proof}

For a uniform hypergraph $\mathcal{H}$, the following is the expression of $C(j)$ by the number of arborescences of the multi-digraphs associated with  $\mathcal{H}$.

%A closed walk $w$ in a multi-digraph $G$ is called an \textit{Eulerian tour} if $w$ traverses each arc of the multi-digraph exactly once.
\begin{thm}\label{waa1}
Let $\mathcal{H}=(V(\mathcal{H}), E(\mathcal{H}))$ be an $m$-uniform hypergraph with $n$ vertices. Then
\begin{equation*}\label{py1}
C(j)=\sum \limits_{d=0}^{\infty}\frac{\mu_{d}(j)}{d!},
\end{equation*}
where
\begin{equation*}
\mu_{d}(j)=(m-1)^{n-1}\sum\limits_{F\in\mathcal{F}_{d}^{(j)}(\mathcal{H})}\frac{t(F)}{\prod\limits_{v\in V(\widetilde{E}(F))/\{j\}}d^{+}(v)},
\end{equation*}
$j\in V(\mathcal{H})$, $d^{+}(v)$ is the outdegree of a vertex $v$ in $D(F)$ and $t(F)$ is the number of arborescences of $D(F)$.
%$E(F)$ is the number of $F'\in F_{d}$ satisfying $E(F')=E(F)$.
\end{thm}
\begin{proof}
Let  $V(\mathcal{H})=\{1,\ldots,n\}$ and $F\in\mathcal{F}_{d}^{(j)}(\mathcal{H})$. Since $F\in\mathcal{F}_{d}^{(j)}(\mathcal{H})$, $D(F)$ is Eulerian and contains the vertex $j$. In order to get the Theorem \ref{waa1}, we first give the representation of $|\textbf{W}(F,j)|$.
Let $\epsilon(F,e_{1})$ be the set of all Eulerian tours starting with a fixed arc $e_{1}$ in $D(F)$. Farrell and Levine \cite{matthew2016multi} proved
\begin{equation}\label{waaa1}
|\epsilon(F,e_{1})|=\frac{\sum\limits_{e\in \widetilde{E}(F)}|\epsilon(F,e)|}{\sum\limits_{v\in V(\widetilde{E}(F))}d^{+}(v)}.
\end{equation}

By Equation (\ref{waaa1}), we have $|\epsilon(F,e_{1})|=|\epsilon(F,e_{2})|, \text{for all} ~e_{1},e_{2}\in \widetilde{E}(F)$. Thus, $$|\epsilon(F,j)|=d^{+}(j)|\epsilon(F,e_{1})|.$$

Let $\mathfrak{E}(F)$ be the set of all Euler circuits in $D(F)$.
There are $|\widetilde{E}(F)||\mathfrak{E}(F)|$ Euler tours in $\widetilde{E}(F)$ \cite{clark2021harary}.
%as we may distinguish any arc from $D(F)$ as the first arc of an Euler circuit.
Then $\sum_{e\in \widetilde{E}(F)}|\epsilon(F,e)|=|\widetilde{E}(F)||\mathfrak{E}(F)|$.

Thus,
$$
|\textbf{W}(F,j)|=\frac{|\epsilon(F,j)|}{b(F)}=\frac{d^{+}(j)\frac{|\widetilde{E}(F)||\mathfrak{E}(F)|}{\sum\limits_{v\in V(\widetilde{E}(F))}d^{+}(v)}}{b(F)}=\frac{d^{+}(j)|\mathfrak{E}(F)|}{b(F)},
$$
where $b(F)$ is the product of the factorials of the multiplicities of all the arcs in $D(F)$.

Tutte et al.\cite{tutte1941on} and Aardenne-Ehrenfest et al.\cite{Aardenne1987Circuits}  proved the BEST Theorem $$|\mathfrak{E}(F)|=t(F)\prod_{v\in V(\widetilde{E}(F))}(d^{+}(v)-1)!.$$
Then
\begin{equation}\label{waa2}
|\textbf{W}(F,j)|=\frac{d^{+}(j)t(F)\prod_{v\in V(\widetilde{E}(F))}(d^{+}(v)-1)!}{b(F)}.
\end{equation}
By Theorem $\ref{z2}$, Equation (\ref{waa2}) and  $c(F)=\prod_{v\in V(\widetilde{E}(F))}d^{+}(v)!$, we get
\begin{align}\label{wei2}
\mu_{d}(j)&=(m-1)^{n-1}\sum\limits_{F\in\mathcal{F}_{d}^{(j)}(\mathcal{H})}\frac{b(F)d^{+}(j)t(F)\prod_{v\in V(\widetilde{E}(F))}(d^{+}(v)-1)!}{c(F)b(F)}\notag\\
&=(m-1)^{n-1}\sum\limits_{F\in\mathcal{F}_{d}^{(j)}(\mathcal{H})}\frac{t(F)}{\prod\limits_{v\in V(\widetilde{E}(F))/\{j\}}d^{+}(v)}.
\end{align}
Substituting Equation (\ref{wei2}) into Equation (\ref{py11}), we obtain the expression of $C(j)$ by the number of arborescences of the multi-digraphs associated with  $\mathcal{H}$.
\end{proof}
Next, we give the following explicit formula of $\mu_{d}(j)$ $(d=0,1,2,\ldots,m)$ for $m$-uniform hypergraphs.
\begin{thm}\label{wa2}
Let $\mathcal{H}=(V(\mathcal{H}),E(\mathcal{H}))$ be an $m$-uniform hypergraph with $n$ vertices. Then
\begin{equation*}
\mu_{d}(j)=\begin{cases}
m^{m-2}(m-1)^{n-m}d(j), &\text{if}~ d=m,\\
0,& \text{if}~ d=1,2,\ldots,m-1,\notag\\
(m-1)^{n-1},& \text{if}~ d=0,
\end{cases}
\end{equation*}
%$$\mu_{m}(j)=m^{m-2}(m-1)^{n-m}d(j),$$
%$$Tr_{3}(\mathcal{A}_\mathcal{H})=\sum \limits_{j=1}^{n}\mu_{3}(j)=\sum \limits_{j=1}^{n}3\cdot 2^{n-3}p_{j}$$
where $j\in V(\mathcal{H})$ and $d(j)$ is the degree of $j$ in $\mathcal{H}$.
\end{thm}
\begin{proof}
 Let $V(\mathcal{H})=\{1,\ldots,n\}$ and $E_{j}(\mathcal{H})=\{e\in E(\mathcal{H})|j\in e\}$. Since $\mathcal{F}_{m}^{(j)}(\mathcal{H})=\{(i_{1}i_{2}\cdots i_{m},i_{2}i_{1}i_{3}i_{4}\cdots i_{m},\ldots,i_{m}i_{1}i_{2}\cdots i_{m-1})|i_{1}i_{2}\cdots i_{m}\in E_{j}(\mathcal{H}), i_{1}i_{2}\cdots i_{m}~\text{is}\\ \text{ascending order} \}$, we have $|\mathcal{F}_{m}^{(j)}(\mathcal{H})|=d(j)$.

For $F\in\mathcal{F}_{m}^{(j)}(\mathcal{H})$, $D(F)$ is a complete digraph with $m$ vertices. We have $\prod_{v\in V(\widetilde{E}(F))/\{j\}}d^{+}(v)=(m-1)^{m-1}$. The number of arborescences of $D(F)$ is $m^{m-2}$ ( Cayley's formula \cite{cayley1889A}).
By Theorem \ref{waa1} and the above discussion, we get $$\mu_{m}(j)=\frac{(m-1)^{n-1}d(j)m^{m-2}}{(m-1)^{m-1}}=m^{m-2}(m-1)^{n-m}d(j).$$

From Lemma \ref{L2} (2), we know $Tr_{d}(\mathcal{A}_\mathcal{H})=0$, $d=1,2,\ldots,m-1$, by Equation (\ref{py}), we know $\sum_{j=1}^n\mu_{d}(j)=Tr_{d}(\mathcal{A}_\mathcal{H})$ and by Theorem \ref{z2}, we get $\mu_{d}(j)\geq0, d=1,2,\ldots,m-1$. Thus, we have $$\mu_{d}(j)=0, d=1,2,\ldots,m-1.$$

By Equation (\ref{A1}), we easily get $\mu_{0}(j)=(m-1)^{n-1}.$
%For $F\in\mathcal{F}_{d}(\mathcal{H}), d=1,2,...,m-1$, $D(F)$ is not Eulerian.  Then, $\mathcal{F}_{d}^{(j)}(\mathcal{H})=\emptyset, d=1,2,...,m-1$. Thus, we have $$\mu_{d}(j)=0, d=1,2,...,m-1.$$
\end{proof}

For an $m$-uniform hypergraph $\mathcal{H}$ with $n$ vertices and $m\geq2$, the subgraph centrality $C(j)=\sum_{d=0}^{\infty}\frac{\mu_{d}(j)}{d!}$, by Theorem \ref{wa2}, we have $\sum_{d=0}^{m}\frac{\mu_{d}(j)}{d!}=(m-1)^{n-1}+\frac{m^{m-3}(m-1)^{n-m-1}d(j)}{(m-2)!},~j=1,2,\ldots,n$.  The ranking from taking $\sum_{d=0}^{m}\frac{\mu_{d}(j)}{d!}$ as centrality is the same as degree centrality for the $m$-uniform hypergraph.

\vspace{3mm}

\noindent
\textbf{Acknowledgements}
\vspace{3mm}
\noindent

%The authors would like to thank the reviewers for giving valuable suggestions.
This work is supported by the National Natural Science Foundation of China (No.
11801115, No. 12071097, No. 12042103 and No. 12242105), the Natural Science Foundation of the
Heilongjiang Province (No. QC2018002) and the Fundamental Research Funds for
the Central Universities.

\section*{References}
\bibliographystyle{unsrt}
\bibliography{spbib}

\begin{thebibliography}{10}

\bibitem{ESTRADA2000713}
E.~Estrada.
\newblock Characterization of 3{D} molecular structure.
\newblock {\em Chem. Phys. Lett.}, 319(5):713--718, 2000.

\bibitem{DELAPENA200770}
J.A. de~la Pe$\tilde{n}$a, I.~Gutman, and J.~Rada.
\newblock Estimating the {E}strada index.
\newblock {\em Linear Algebra Appl.}, 427(1):70--76, 2007.

\bibitem{Gutman2008Lower}
I.~Gutman.
\newblock Lower bounds for {E}strada index.
\newblock {\em Publ. Inst. Math.}, 83:1--7, 2008.

\bibitem{2008On}
B.~Zhou.
\newblock On {E}strada index.
\newblock {\em MATCH Commun. Math. Comput. Chem.}, 60(2):485--492, 2008.

\bibitem{2011The}
Z.~Du and B.~Zhou.
\newblock The {E}strada index of trees.
\newblock {\em Linear Algebra Appl.}, 435(10):2462--2467, 2011.

\bibitem{2012Estrada}
Z.~Chen, Y.~Fan, and W.~Du.
\newblock {E}strada index of random graphs.
\newblock {\em MATCH Commun. Math. Comput. Chem.}, 68(3):815--823, 2012.

\bibitem{molecules}
E.~Estrada, J.A. Rodr\'{i}guez-Vel\'{a}zquez, and M.~Randi\'{c}.
\newblock Atomic branching in molecules.
\newblock {\em Int. J. Quantum Chem.}, 106(4):823--832, 2006.

\bibitem{PhysRevE.71.056103}
E.~Estrada and J.A. Rodr\'{\i}guez-Vel\'azquez.
\newblock Subgraph centrality in complex networks.
\newblock {\em Phys. Rev. E}, 71:056103, 2005.

\bibitem{PhysRevE.72.046105}
E.~Estrada and J.A. Rodr\'{\i}guez-Vel\'azquez.
\newblock Spectral measures of bipartivity in complex networks.
\newblock {\em Phys. Rev. E}, 72:046105, 2005.

\bibitem{1980Spectra}
D.~Cvetkovi{\'c}, M.~Doob, and H.~Sachs.
\newblock {\em Spectra of Graphs}.
\newblock Academic {P}ress, 1980.

\bibitem{2011Thee}
E.~Estrada.
\newblock {\em The {S}tructure of {C}omplex {N}etworks: {T}heory and
  {A}pplications}.
\newblock Oxford {U}niversity {P}ress, 2011.

\bibitem{qi2005eigenvalues}
L.~Qi.
\newblock Eigenvalues of a real supersymmetric tensor.
\newblock {\em J. Symbolic Comput.}, 40(6):1302--1324, 2005.

\bibitem{lim2005singular}
L.~Lim.
\newblock Singular values and eigenvalues of tensors: a variational approach.
\newblock In {\em 1st IEEE International Workshop on Computational Advances in
  Multi-Sensor Adaptive Processing}, pages 129--132. IEEE, 2005.

\bibitem{cooper2012spectra}
J.~Cooper and A.~Dutle.
\newblock Spectra of uniform hypergraphs.
\newblock {\em Linear Algebra Appl.}, 436(9):3268--3292, 2012.

\bibitem{cooper2}
Y.~Fan, T.~Huang, Y.~Bao, C.~Zhuan-Sun, and Y.~Li.
\newblock The spectral symmetry of weakly irreducible nonnegative tensors and
  connected hypergraphs.
\newblock {\em Trans. Amer. Math. Soc.}, 372(3):2213--2233, 2019.

\bibitem{doi:10.1137/21M1404740}
G.~Gao, A.~Chang, and Y.~Hou.
\newblock Spectral {R}adius on {L}inear $r$-{G}raphs without {E}xpanded
  ${K}_{r+1}$.
\newblock {\em SIAM J. Discrete Math.}, 36(2):1000--1011, 2022.

\bibitem{2017Tensor}
L.~Qi and Z.~Luo.
\newblock {\em Tensor {A}nalysis: {S}pectral {T}heory and {S}pecial {T}ensors}.
\newblock S{I}{A}{M}, 2017.

\bibitem{clark2021harary}
G.~Clark and J.~Cooper.
\newblock A {Harary-Sachs} theorem for hypergraphs.
\newblock {\em J. Combin. Theory Ser. B}, 149:1--15, 2021.

\bibitem{2011Analogue}
A.~Morozov and S.~Shakirov.
\newblock Analogue of the identity {Log Det = Trace Log} for resultants.
\newblock {\em J. Geom. Phys.}, 61(3):708--726, 2011.

\bibitem{hu2013determinants}
S.~Hu, Z.~Huang, C.~Ling, and L.~Qi.
\newblock On determinants and eigenvalue theory of tensors.
\newblock {\em J. Symbolic Comput.}, 50:508--531, 2013.

\bibitem{shao2015some}
J.~Shao, L.~Qi, and S.~Hu.
\newblock Some new trace formulas of tensors with applications in spectral
  hypergraph theory.
\newblock {\em Linear Multilinear Algebra}, 63(5):971--992, 2015.

\bibitem{doi:10.1080/03081087.2021.1953431}
L.~Chen, C.~Bu, and J.~Zhou.
\newblock Spectral moments of hypertrees and their applications.
\newblock {\em Linear Multilinear Algebra}, 70(21):6297--6311, 2022.

\bibitem{FAN2021112329}
Y.~Fan, M.~Li, and Y.~Wang.
\newblock The cyclic index of adjacency tensor of generalized power
  hypergraphs.
\newblock {\em Discrete Math.}, 344(5):112329, 2021.

\bibitem{2017On}
G.~Clark and J.~Cooper.
\newblock On the adjacency spectra of hypertrees.
\newblock {\em Electron. J. Combin.}, 25(2), 2018.
\newblock Article ID \#P2.48.

\bibitem{Hu2013Cored}
S.~Hu, L.~Qi, and J.~Shao.
\newblock Cored hypergraphs, power hypergraphs and their {L}aplacian
  {H}-eigenvalues.
\newblock {\em Linear Algebra Appl.}, 439:2980--2998, 2013.

\bibitem{1978The}
I.~Gutman.
\newblock The {E}nergy of a graph.
\newblock {\em Ber. Math-Statist. Sekt. Forschungszentrum Graz}, 103:1--22,
  1978.

\bibitem{BAMDAD2010739}
H.~Bamdad, F.~Ashraf, and I.~Gutman.
\newblock Lower bounds for {E}strada index and {L}aplacian {E}strada index.
\newblock {\em Appl. Math. Lett.}, 23(7):739--742, 2010.

\bibitem{matthew2016multi}
M.~Farrell and L.~Levine.
\newblock Multi-{E}ulerian tours of directed graphs.
\newblock {\em Electron. J. Combin.}, 23(2), 2016.
\newblock Article ID \#P2.21.

\bibitem{tutte1941on}
W.T. Tutte and C.A.B. Smith.
\newblock On unicursal paths in a network of degree 4.
\newblock {\em Amer. Math. Monthly}, 48(4):233--237, 1941.

\bibitem{Aardenne1987Circuits}
T.~van Aardenne-Ehrenfest and N.G. de~Bruijin.
\newblock Circuits and {T}rees in {O}riented {L}inear {G}raphs.
\newblock Master's thesis, Birkh$\rm\ddot{a}$user Boston, 1987.

\bibitem{cayley1889A}
A.~Cayley.
\newblock A theorem on trees.
\newblock {\em Q. J. Math.}, 23:376--378, 1889.

\end{thebibliography}
\end{spacing}
\end{document}